\documentclass{amsart}
\usepackage[cp1251]{inputenc}
\usepackage[english]{babel}
\usepackage{amsmath}
\usepackage{amssymb}
\usepackage{amsfonts}

\newtheorem{theorem}{Theorem}

\newtheorem{definition}{Definition}

\begin {document}

\thispagestyle{empty}

\title {On maximal extensions of nilpotent Lie algebras}

\author {V.V.Gorbatsevich}

\maketitle

The article studies finite-dimensional Lie algebras $L$ with a given nilradical. Lie algebras having a given (and considered up to isomorphism) nilradical $N$ will be called extensions of the nilpotent Lie algebra $N$. If Lie algebra solvable then the extension will be called solvable.

Such extensions can be useful in various fields of mathematics and physics. In particular, the description of solvable Lie groups and algebras is useful in the study of symmetries of differential and other equations, as well as in the study of representations of solvable Lie groups. In physics, solvable Lie algebras and groups are used in string theory and other theories of elementary particles, as well as in some higher-dimensional theories in cosmology.

We will assume that the field $k$ of the Lie algebras we are considering is of characteristic 0 and sometimes is algebraically closed. For Lie algebras over a field $k$ we need the notion of a toric subalgebra Lie, or, in other words, an Abelian Lie subalgebra consisting of semisimple (i.e., completely reducible in the adjoint representation) elements. As is known, for an algebraically closed field of characteristic 0, all maximal toric Lie subalgebras are conjugate. For non-closed fields this is no longer the case --- for example, in the simple Lie algebra $sl_2(\bf R)$ there are two nonconjugate (one-dimensional) toric Lie subalgebras --- compact (isomorphic to $so(2)$) and split (isomorphic Lie subalgebra of diagonal matrices with trace 0).

The general information about Lie algebras used in the article can be found, for example, in \cite {VGO}, and properties of algebraic Lie algebras, for example, in \cite {H}. For a Lie algebra $L$, by $Der(L)$ we denote the Lie algebra of its derivations, and by $IDer(L)$ the Lie algebra of its inner derivations (induced by the action of the Lie algebra itself on itself by means of the adjoint representation $ad_L$).

Let $N$ be some nilpotent Lie algebra over the field $k$. Consider all possible finite-dimensional Lie algebras $L$ over the field $k$ whose nilradical is isomorphic to $N$. It is clear that such Lie algebras $L$ exist and there are infinitely many of them. For example, such are all Lie algebras of the form $N \oplus S$ --- direct sums of the Lie algebra $N$ and an arbitrary semisimple Lie algebra $S$. Note that in this example the adjoint action of the Lie algebra $S$ on the nilradical $N$ is trivial. In general, the direct sum of an arbitrary Lie algebra and any semisimple Lie algebra does not change the nilradical of the original Lie algebra.

Let us indicate one more way of constructing Lie algebras, moreover, solvable ones, with a given nilradical $N$. Consider the Lie algebra $Der(N)$ of the derivation of the Lie algebra $N$. The Lie algebra $Der(N)$ is known to be algebraic (the nilpotency of the Lie algebra $N$ is not essential here). Like any algebraic Lie algebra, the Lie algebra $Der(N)$ has a Chevalley decomposition $Der(N)=P+U$ into a semidirect sum of a maximal reductive Lie subalgebra $P$ and a nilradical $U$ (consisting of nilpotent elements). The reductive Lie subalgebra (for an algebraically closed field $k$ it is unique up to conjugation in $Der(N)$) decomposes --- and uniquely --- into a direct sum $P=S \oplus T$ of a semisimple Lie algebra $S$ and a central Lie subalgebra $T$ (toric, that is, consisting of semisimple elements). If $T_S$ is some Cartan Lie subalgebra in $S$ (or, which is the same here, a maximal Abelian Lie subalgebra consisting of semisimple elements), then $T_S \oplus T$ is a maximal toric subalgebra in the Lie algebra $ Der(N)$ (maximal Abelian Lie subalgebra in $Der(N)$ consisting of semisimple elements). Consider the Lie algebra $L= (T_S \oplus T) + N$, the semidirect sum of the Lie subalgebra $T_S \oplus T$ and the ideal $N$, corresponding to the natural action of the Lie subalgebra $T_S \oplus T \subset Der(N) $ on $N$ by derivations. It is clear that the Lie algebra $L$ is solvable and its nilradical is isomorphic to the original nilpotent Lie algebra $N$. Similarly, considering arbitrary reductive Lie subalgebras in $Der(N)$ instead of $T_S \oplus T$ and forming their semidirect sums with $N$, we generally obtain nonsolvable Lie algebras with a nilradical isomorphic to $N$.

Of interest (from a negative point of view) is the case of characteristically nilpotent Lie algebras (that is, those whose Lie algebra of their derivations is nilpotent). It is easy to see that such Lie algebras do not have "non-trivial" (that is, not isomorphic to the original Lie algebras) solvable extensions.

In this article, we are interested in the maximal ones among Lie algebras $L$ having some fixed nilradical $N$. Solvable Lie algebras $R$ with a given nilradical are often considered as a special, but very important, particular case. Since there are many such Lie algebras $L$ and $R$, we select the maximal ones among them. At the same time, maximality can be understood in two different senses. The first is the maximality with respect to inclusions, i.e. a maximal Lie algebra is not isomorphic to a proper subalgebra in any other Lie algebra with a similar property (in our case, with the same nilradical). The second is maximality in dimension (when among all Lie algebras with a given property, those that have the largest dimension are singled out). Moreover, if we consider arbitrary (not only solvable) finite-dimensional Lie algebras with a given nilradical, then it is clear that among such Lie algebras there are no maximal ones in any of the indicated senses. This is due to the fact that adding an arbitrary semisimple direct summand to the Lie algebra $L$ does not change the nilradical. Therefore, here it is necessary to single out a special class of Lie algebras with a given nilradical --- which, by analogy with the terminology introduced in \cite {HO}, will - and as extensions - be called exact. Namely, a Lie algebra $L$ with nilradical $N$ is said to be exact if the kernel of the adjoint action of the semisimple part $S$ of the Lie algebra $L$ on $N$ is trivial (which is known to be equivalent to the triviality of the kernel of the action of $S$ on radical $R$ of the Lie algebra $L$). This condition is easily understood to be equivalent to the fact that the Lie algebra $L$ cannot be represented as a direct sum of a semisimple Lie algebra and some other Lie algebra (whose nilradical is obviously isomorphic to $N$). If, however, we consider only solvable Lie algebras with a given nilradical, then such additional conditions are no longer needed.

For an arbitrarily taken property, these two notions of maximality are, generally speaking, different.

In this paper, we will use the notion of Lie algebra splitting, in particular, the Malcev splitting introduced in \cite {Ma}. It is described in a slightly different form in \cite {G1}. Here it should be noted that the fundamental work of A.I.Malcev \cite {Ma} turned out to be almost completely undeservedly forgotten. There are very few references to it. It is usually mentioned only in connection with the fact that there the classification of complex solvable Lie algebras is reduced to the classification of complex nilpotent Lie algebras. This statement is not entirely accurate, because, for such a reduction, as shown in that article, it is also necessary to describe all possible orbits of some linear actions of some nilpotent Lie groups, which itself is a non-trivial problem. But, on the other hand, many important results on solvable (and not only on solvable) Lie algebras were obtained in this paper.

\begin{definition} A Lie algebra $L$ is said to be split if it admits a decomposition of the form $L = S + T + U$, where $U$ is its nilradical, $S$ is its semisimple part, the Lie subalgebra $T$ is Abelian, and $T$ centralizes the Lie subalgebra $S$, and for any $X \in T$ the linear operator $ad_L(X)$ is nontrivial and semisimple.

An embedding of a Lie algebra into a split Lie algebra, as well as this split Lie algebra itself, we will call a splitting of $L$.
\end{definition}

It follows from the definition of the splitting of the Lie algebra $L$ that $T \cap U = \{0 \}$ (because elements from $T\cap U$ act on $L$ semisimply and nilpotently, and therefore trivially), and also that $T$ is isomorphic to a subalgebra in $Der (U)$ --- the Lie algebra of derivations of the Lie algebra $U$.

Note that all algebraic Lie algebras are split (this follows from the Chevalley decomposition --- see, for example, \cite {H}).

\begin{definition} Let $L$ be a Lie algebra. A Malcev splitting is an embedding $\alpha: L \hookrightarrow M(L)$ into a split algebra $M(L) =S+T+U$ such that $M(L)$ is the semidirect sum of the subalgebra $T$ and the ideal $\alpha(L)$ and $\alpha(L) +U=M(L)$.
\end{definition}

The Malcev's article also contains an additional minimality condition for a split Lie algebra containing the original Lie algebra. We will call the corresponding splitting as Malcev-type splitting.

In \cite {G1} it is proved that for an arbitrary finite-dimensional Lie algebra $L$ a splitting exists over a field of characteristic 0 (for complex Lie algebras this was first proved by A.I. Malcev \cite {Ma}) and it is unique (with respect to the naturally defined splitting isomorphism). Here we briefly present the construction of the Malcev splitting for Lie algebras (it is given in more detail --- for the parallel notion of the splitting of Lie groups --- in \cite {G1}).

Let $L=S+R$ be the Levi decomposition of the Lie algebra $L$. Consider an adjoint representation $ad_L : L \to gl(L)$ and set $L^\ast =ad_L(L)$. Further, let $<L^\ast>$ be the algebraic closure of the Lie subalgebra $L^\ast$ in $gl(L)$ (that is, the smallest algebraic Lie subalgebra containing $L^\ast$). Since the Lie algebra $<L^\ast>$ is algebraic, we have the Chevalley decomposition for it: $<L^\ast> = S^\ast +T^\ast+U^\ast$, where $U^\ast $ is a nilpotent radical, the Lie subalgebra $S^\ast$ is semisimple (Levy factor), and the Lie subalgebra $T^\ast$ is Abelian and consists of semisimple (i.e., completely reducible over $k$) elements.
Moreover, the Lie subalgebra $<L^\ast>$ is contained in $Der(L)$, since the Lie subalgebra $Der(L)$ is obviously algebraic, while $L^\ast$ is contained in it.
Consider the semidirect sum $W^\ast=S^\ast +U^\ast$, then $<L^\ast> = T^\ast+W^\ast$, and $T^\ast \cap W^ \ast = \{ 0 \}$. Let $f : T^\ast+W^\ast \to T^\ast$ be the natural epimorphism with kernel $W^\ast$. Denote $T=f(L)$, then it is clear that $T \subset Der(L)$ and therefore we can form a semidirect sum $T+L$ corresponding to the natural action of $T$ on $L$. The resulting Lie algebra $T+L$ (with the natural embedding of $L$ into $M(L)$) we call the Malcev splitting for $L$ and denote it by $M(L)$. In \cite {G1} it is proved that this construction gives a some kind of splitting. More precisely, a similar assertion is proved there for arbitrary connected Lie groups, and our assertion for Lie algebras is equivalent to the case when these Lie groups are simply connected. Note that the Malcev splitting for a given Lie algebra is unique up to a naturally defined notion of isomorphism of splittings \cite {G1}.

In \cite {Sn} it was proved that $\dim R/N$ does not exceed $\dim{N/[N,N]}$ (a slightly different terminology is used in that paper). We will now prove a generalization of this assertion.

Note that the Lie algebra $L/N$ is always reductive. Its toric rank (the dimension of the maximal toric Lie subalgebra) will be denoted by $r_t(L/N)$.

\begin{theorem} Let $L$ be an arbitrary Lie algebra over an algebraically closed field of characteristic 0, witch is exact, and $N$ be its nilradical.
Then $r_t(L/N) \le \dim N/[N.N]$

If $L$ is solvable, then the exactness condition is, of course, redundant. Here the rank of the toric (Abelian) Lie algebra $L/N$ is equal to its dimension.
\end{theorem}

\begin{proof}
Consider the Malcev splitting for the Lie algebra $L=S+R$. We have $M(L)=T+L$, where $T$ is a toric Lie subalgebra in $Der(L)$.

In what follows, we will use the following well-known fact: for a nilpotent Lie algebra $N$, the minimum number of generators is equal to $\dim(N/[N,N])$. The same can be expressed as follows: if $V$ is a complementary subspace of $[N,N]$ in $N$, then $V$ generates a Lie algebra $N$.

The adjoint action of the reductive subalgebra $S\oplus T \subset M(L)$ on $L$ preserves $N$ (as the characteristic ideal in $L$) and the ideal $[N,N]$. Moreover, the action of the subalgebra $S \oplus T$ on $N$ generates --- due to its reductivity --- the action on some subspace $V \subset N$ complementary to $[N,N]$ in $N$. If for some nonzero $X \in S \oplus T$ it turns out that its action on $V$ is trivial, then its action on the entire Lie algebra $N$ will also be trivial (since $V$ generates $N$ ). On the other hand, the set of all those elements from $S\oplus T \subset M(L)$ whose action on $V$ is trivial is an ideal in $S\oplus T \subset M(L)$. But ideals are in $T\oplus S$, obviously. have the form $S^\prime \oplus T^\prime$, where $S^\prime$ is an ideal in $S$ and $T^\prime$ is some subalgebra of $S$.

Assume that the ideal $S^\prime$ for the kernel of the indicated action is non-trivial. But then the fact that $V$ generates the Lie algebra $N$ implies that the Lie algebra $L$ is not exact, which contradicts our assumption. Therefore $S^\prime= \{0\}$ and the kernel of the action lies in $T$. But due to the Malcev splitting construction, this can only be true if $T^\prime = \{0\}$. This proves that the homomorphism of the Lie algebra $T \oplus S$ to $gl(V)$ is exact, i.e., it has the trivial kernel. But then it is obvious that $r_t(L/N) \le \dim (N/[N,N])$.

For solvable $R$, the assertion of the theorem immediately follows from what has already been proved above.
\end{proof}

We pass to the consideration of maximal extensions of nilpotent Lie algebras. Let's start with the case of solvable extensions.
L. $\check{S}$noble in \cite {Sn} (see also \cite {SnW}) put forward a conjecture about the form of maximal (in dimension) solvable extensions of complex nilpotent Lie algebras $N$. It was assumed that all of them are semidirect sums of maximal Abelian Lie subalgebras $T\subset Der(N)$ consisting of semisimple elements and the initial Lie algebra $N$. In particular, it was assumed that for every nilpotent Lie algebra its maximal (in dimension) extension is unique up to isomorphism. The validity of this hypothesis was confirmed by the fact that it was indeed tested in a number of special cases (a list of relevant articles can be found in \cite {Sn}). For example, it is true for all nilpotent Lie algebras of dimension $\le 6$. It is also true for a number of classes of nilpotent Lie algebras of arbitrary dimension. However, as will be shown below, this conjecture is not true in the general case --- a statement will be given on the basis of which a counterexample will be given --- a nilpotent Lie algebra of dimension 8, for which two non-isomorphic maximal solvable extensions will be indicated.

Note that the analogue of $\check{S}$noble's conjecture from \cite {Sn} for fields that are not algebraically closed is false for a very simple reason. The point is that for such fields the maximal toric Lie subalgebras are very often not conjugate to each other. This makes it easy to construct counterexamples to Snoble's conjecture. For example, for $N={\bf R}^2$ there are two maximal (with respect to embedding) solvable extensions (they correspond to two toric Lie subalgebras). These are $(so(2)\oplus C)+{\bf R}^2$ (where $C$ is a subalgebra of scalar matrices) and the Lie algebra ${\bf R}^2+{\bf R}^ 2$ corresponding to the two-dimensional Abelian split Lie subalgebra in $gl_2(\bf R)$ formed by all diagonal matrices.

We now turn to constructing a counterexample to $\check{S}$noble's conjecture formulated above.

Let $N_7$ be a seven-dimensional nilpotent Lie algebra for which the Lie algebra $Der(N)$ is also nilpotent (such Lie algebras are called characteristically nilpotent; as is well known, they exist only starting from dimension 7). Consider a nilpotent Lie algebra $N = k\oplus N_7$, where $k$ is a one-dimensional Abelian nilpotent Lie algebra. It is clear that the Lie algebra $Der(k)$ is isomorphic to the one-dimensional toric Lie algebra $gl_1(k)$.

We need to describe the Lie algebra $Der(N)$. To do this, we use the description of Lie algebras of derivations of direct sums of Lie algebras obtained in \cite {T}. In our case, $N$ has two direct summands $k$ and $N_7$. Therefore, due to \cite {T}, we have a representation of $Der(N)$ as a sum of four subspaces:

$Der(N)= (Der(k) \oplus Der(N_7))+ (Der(k, N_7)+Der(N_7,k))$

where $Der(k,N_7)$ is the set of linear mappings $k \to Z(N_7)$ to the center $Z(N_7)$ (nontrivial due to the nilpotency of the Lie algebra $N_7$), and $Der(N_7 ,k)$ is the set of linear mappings $N_7/[N_7,N_7] \to k$. Here we have specified the general statement from \cite {T} as applied to the case under consideration.

Note that the maximal reductive Lie subalgebra in $Der(N)$ is the direct sum of the maximal reductive subalgebras in the algebras of derivations of the Lie algebras $k$ and $N_7$ (obviously, algebraic, in particular --- split). This fact was noted in \cite {HO}.

We are interested in the Lie subalgebra $Der(k) \oplus Der (N_7)$ in $Der(N)$. The Lie algebra $Der(k)$ is a toric Lie algebra (we denote it by $T$ and let $X \in T$ be a nonzero vector), and in $Der(N_7)$ we choose and fix some outer derivation $d$ (as noted above, such an outer derivation always exists for nilpotent Lie algebras). For concreteness, we take the central derivation, which takes the element $X_1$ to the generator $X_7$ of the (one-dimensional) center, and the other elements of the basis to 0. It is clear that the subalgebra $T$ permutes with $d$. Consider two solvable extensions of the nilpotent Lie algebra $N$: $R_1=r_2\oplus N_7$ (it is obtained by deriving $X$), and $R_2$ is the extension corresponding to deriving $X+d$. Obviously, these extensions are maximal. And these two Lie algebras will not be isomorphic with each other. To see this, we note that the Lie algebra $R_1$ is split, and therefore its Maltsev splitting coincides with itself (and has dimension 9). As regards the Lie algebra $R_2$, it is easy to verify that its Malcev-type splitting has dimension 10. Thus, the Malcev-type splittings of the Lie algebras $R_1,R_2$ have different dimensions, and, due to the uniqueness of the Malcev-type splitting, they therefore do not will be isomorphic. We get that the Lie algebra $N=k \oplus N_7$ has two maximal solvable non-isomorphic extensions.

But there is another (purely computational) way of proving the non-isomorphism of $R_1$ and $R_2$, which was reported to the author by B.Omirov. If we calculate the Lie algebras of derivations of the Lie algebras $R_1,R_2$, for example, for the nilpotent Lie algebra of G. Favre, then it turns out that the dimensions of these derivation algebras are different --- they are equal to 13 and 12, respectively. For calculations, one can use mathematical packages containing  symbolic calculations (the author used the GAP package).

The constructed example shows that the above conjecture of L. Shnobl turned out to be incorrect in the general case. We note that completely similar examples are constructed for some of those nilpotent Lie algebras that can be represented as a direct sum of some non-trivial nilpotent Lie algebra (for which we took the one-dimensional one above) and some characteristically nilpotent one (which exist in any dimension, starting from 7) .

In what follows, some information  about arbitrary exact extensions of nilpotent Lie algebras is given . Recall that there are no maximal extensions of nilpotent Lie algebras on which the exactness condition is not imposed --- namely, there are extensions of an arbitrary nilpotent Lie algebra in which it is a nilradical and whose dimensions have no upper bound.

Let $L=S+R$ be some exact extension of a nilpotent Lie algebra $N$. Consider the Malcev splitting $M(L)=S+T+U$. Here the semisimple part of the Lie algebra $L$ is denoted by the same letter $S$ as the semisimple part of the Lie algebra $M(L)$. The point is that, as shown in \cite {Ma},\cite {G1}, these two semisimple parts are isomorphic to each other. Maximal exact extensions of the Lie algebra $N$ will obviously have the form $S+T+N$ ($S+T$ is the maximal reductive Lie subalgebra in $Der(N)$). We will call them standard extensions of the nilpotent Lie algebra $N$. As for solvable extensions, in addition to these standard extensions, in some cases there may also be maximal non-standard exact extensions.

\end{document}